\documentclass[12pt,italian,english]{amsart}
\usepackage[T1]{fontenc}
\usepackage[latin9]{inputenc}
\usepackage{geometry}
\usepackage{babel}
\usepackage{amstext}
\usepackage{amsthm}
\usepackage{amssymb}
\usepackage[unicode=true,pdfusetitle,
 bookmarks=true,bookmarksnumbered=false,bookmarksopen=false,
 breaklinks=false,pdfborder={0 0 1},backref=false,colorlinks=false]
 {hyperref}

\makeatletter
\theoremstyle{plain}
\newtheorem{thm}{\protect\theoremname}

\makeatother

  \addto\captionsenglish{\renewcommand{\theoremname}{Theorem}}
  \addto\captionsitalian{\renewcommand{\theoremname}{Teorema}}
\providecommand{\theoremname}{Theorem}

\begin{document}

\title{Some identities involving the Cesàro average of Goldbach numbers}

\author{Marco Cantarini}
\begin{abstract}
Let $\Lambda\left(n\right)$ be the von Mangoldt function and $r_{G}\left(n\right):=\sum_{m_{1}+m_{2}=n}\Lambda\left(m_{1}\right)\Lambda\left(m_{2}\right)$
be the counting function for the numbers that can be written as sum
of two primes (that we will call ``Goldbach numbers'', for brevity)
and let $\widetilde{S}\left(z\right):=\sum_{n\geq1}\Lambda\left(n\right)e^{-nz}$,
with $z\in\mathbb{C}$, $\mathrm{Re}\left(z\right)>0$. In this paper
we will prove the identity
\begin{align*}
\widetilde{S}\left(z\right)= & 2e^{-2z}+\frac{e^{-2z}}{z}-\sum_{\rho}z^{-\rho}\Gamma\left(\rho\right)\\
+ & \sum_{\rho}\left(z^{-\rho}\gamma\left(\rho,2z\right)-\frac{2^{\rho}e^{-z}}{\rho}\right)-\frac{\zeta^{\prime}}{\zeta}\left(0\right)e^{-2z}\\
+ & \frac{e^{-z}}{2}\left(-i\pi\textrm{sgn}\left(\textrm{Im}\left(-z\right)\right)+\textrm{Ei}\left(-z\right)\right)\\
+ & \frac{e^{z}}{2}\left(-i\pi\textrm{sgn}\left(\textrm{Im}\left(-z\right)\right)-\log\left(3\right)e^{-3z}+\textrm{Ei}\left(-3z\right)\right)\\
+ & i\pi\textrm{sgn}\left(\textrm{Im}\left(-z\right)\right)+\log\left(2\right)e^{-2z}-\textrm{Ei}\left(-2z\right)
\end{align*}
where $\gamma\left(\rho,2z\right)$ is the lower incomplete Gamma
function, $\rho=\beta+i\gamma$ runs over the non-trivial zeros of
the Riemann Zeta function and 
\[
\mathrm{Ei}\left(z\right)=-\int_{-z}^{\infty}\frac{e^{-t}}{t}dt,\ \left|\arg\left(z\right)\right|<\pi
\]
 is the Exponential integral function. In addition we will prove that
\begin{align*}
\sum_{n\leq N}r_{G}\left(n\right)\left(N-n\right)= & \frac{N^{3}}{6}-2\sum_{\rho}\frac{\left(N-2\right)^{\rho+2}}{\rho\left(\rho+1\right)\left(\rho+2\right)}\\
+ & \sum_{\rho_{1}}\sum_{\rho_{2}}\frac{\Gamma\left(\rho_{1}\right)\Gamma\left(\rho_{2}\right)}{\Gamma\left(\rho_{1}+\rho_{2}+2\right)}N^{\rho_{1}+\rho_{2}+1}-2\sum_{\rho_{1}}\sum_{\rho_{2}}\frac{N^{\rho_{1}+\rho_{2}+1}}{\rho_{1}\rho_{2}}B_{2/N}\left(\rho_{1}+1,\rho_{2}+1\right)\\
+ & F\left(N\right)
\end{align*}
where $N>4$ is a natural number, $B_{x}\left(a,b\right)$ is the
incomplete Beta function and $F\left(N\right)$ is a sum of (explicitly
calculate) elementary functions, dilogarithms and sums over non-trivial
zeros of the Riemann Zeta function involving the incomplete Beta function.
Furthemore we will prove that
\[
F\left(N\right)=O\left(N^{2}\right)
\]
as $N\rightarrow\infty.$
\end{abstract}

\maketitle
\subjclass 2010 Mathematics Subject Classification:{ Primary 11P32; Secondary 44A10}


\section{Introduction}

In this paper we prove some identities inherent the study to Cesàro
average of Goldbach numbers. The main result in this direction is
the paper of Languasco e Zaccagnini \cite{Langzac}. The authors proved
the formula
\begin{align}
\sum_{n\leq N}r_{G}\left(n\right)\frac{\left(N-n\right)^{k}}{\Gamma\left(k+1\right)}= & \frac{N^{2+k}}{\Gamma\left(k+3\right)}-2\sum_{\rho}\frac{\Gamma\left(\rho\right)}{\Gamma\left(\rho+k+2\right)}N^{\rho+k+1}\nonumber \\
+ & \sum_{\rho_{1}}\sum_{\rho_{2}}\frac{\Gamma\left(\rho_{1}\right)\Gamma\left(\rho_{2}\right)}{\Gamma\left(\rho_{1}+\rho_{2}+k+1\right)}N^{\rho_{1}+\rho_{2}+k}+O_{k}\left(N^{k+1}\right)\label{langzac formula}
\end{align}
for $k>1$, where $N$ is a large natural number and $\rho$, with
or without subscripts, runs over the non-trivial zeros of the Riemann
Zeta function $\zeta\left(s\right)$. The parameter $k$ plays a central
role in the study of these problems and we would like to keep it as
small as possible since, for $k=0$, the Cesàro weight disappears.
The method of Languasco and Zaccagnini allows to study many types
of additive problems; see \cite{Canta1}, \cite{Canta2} and \cite{LangZac2}.
The crucial points of their technique are the relation
\begin{equation}
\sum_{n\leq N}r_{G}\left(n\right)\frac{\left(N-n\right)^{k}}{\Gamma\left(k+1\right)}=\frac{1}{2\pi i}\int_{\left(1/N\right)}e^{Nz}z^{-k-1}\widetilde{S}\left(z\right)^{2}dz\label{identit=0000E0 fondamentale}
\end{equation}
which holds for $k>0,$ where
\[
\widetilde{S}\left(z\right):=\sum_{n\geq1}\Lambda\left(n\right)e^{-nz},\,z=a+iy,\,\mathrm{with}\,a>0,\,y\in\mathbb{R}.
\]
is the power series of the Von Mangoldt function and $\int_{\left(1/N\right)}$
means $\int_{1/N-i\infty}^{1/N+i\infty}$, and the estimation
\begin{equation}
\widetilde{S}\left(z\right)=\frac{1}{z}-\sum_{\rho}z^{-\rho}\Gamma\left(\rho\right)+E\left(a,y\right)\label{eq:asint s tilda zac}
\end{equation}
where
\begin{equation}
E\left(a,y\right)\ll\begin{cases}
1+\left|z\right|^{1/2}, & \left|y\right|\leq a\\
1+\left|z\right|^{1/2}\left(1+\log^{2}\left(\left|y\right|/a\right)\right), & \left|y\right|>a.
\end{cases}\label{errore s tilda}
\end{equation}
The proof of \eqref{eq:asint s tilda zac} is present in \cite{Langzac}.
We recall that the study of $\widetilde{S}\left(z\right)$ is classical,
see for example \cite{HardLit}. One of the Recently Goldston and
Yang \cite{GolYan} proved, assuming RH, that 
\begin{equation}
\sum_{n\leq N}r_{G}\left(n\right)\left(1-\frac{n}{N}\right)=\frac{N^{2}}{6}-2\sum_{\rho}\frac{\Gamma\left(\rho\right)}{\Gamma\left(\rho+3\right)}N^{\rho+1}+O\left(N\right)\label{eq:GoldeYan}
\end{equation}
which correspond to what Languasco and Zaccagnini's formula implies
if one take $k=1$. Languasco and Zaccagnini \cite{LangZacNew} also
proved (assuming RH) that
\begin{align*}
\sum_{n=N-H}^{N+H}r_{G}\left(n\right)\left(1-\frac{\left|n-N\right|}{H}\right)= & HN-\frac{2}{H}\sum_{\rho}\frac{\left(N+H\right)^{\rho+2}-2N^{\rho+2}+\left(N-H\right)^{\rho+2}}{\rho\left(\rho+1\right)\left(\rho+2\right)}\\
+ & O\left(N\log^{2}\left(\frac{2N}{H}\right)+H\log^{2}\left(N\right)\log\left(2H\right)\right)
\end{align*}
where $N\geq2,\,1\leq H\leq N$. $k.$ In a very recent paper Brüdern,
Kaczorowski and Perelli \cite{Perelli} were able to find an explicit
formula which holds for all $k>0.$ Similar averages of arithmetical
functions are common in literature, see, e.g., Chandrasekharan - Narasimhan
\cite{Chana} and Berndt \cite{Ber} who built on earlier classical
work. In In this paper we will prove an explicit formula for $\widetilde{S}\left(z\right)$
which will allow us to find a smaller error respect to \eqref{errore s tilda}.
Kunik and Lucht \cite{Kunik} wrote an explicit formula for $\widetilde{S}\left(z\right)$
but in their closed form there is a term regarding another summation
of Von Mangoldt function, i.e. the series $\sum_{n\geq1}\Lambda\left(n\right)\left(1-\exp\left(-z/n\right)\right)/n$.
We will able to rewrite this series in terms of the lower incomplete
Gamma function and sum over non-trivial zeros of the Riemann Zeta
function. We first prove the following
\begin{thm}
\label{lem:5}Let $z=a+iy,\,a>0,\,y\in\mathbb{R}.$ Let us consider
the function
\[
\widetilde{S}\left(z\right)=\sum_{m\geq1}\Lambda\left(m\right)e^{-mz}.
\]
Then 
\begin{align*}
\widetilde{S}\left(z\right)= & 2e^{-2z}+\frac{e^{-2z}}{z}-\sum_{\rho}z^{-\rho}\Gamma\left(\rho\right)\\
+ & \sum_{\rho}\left(z^{-\rho}\gamma\left(\rho,2z\right)-\frac{2^{\rho}e^{-z}}{\rho}\right)-\frac{\zeta^{\prime}}{\zeta}\left(0\right)e^{-2z}\\
+ & \frac{e^{-z}}{2}\left(-i\pi\mathrm{sgn}\left(\mathrm{Im}\left(-z\right)\right)+\mathrm{Ei}\left(-z\right)\right)\\
+ & \frac{e^{z}}{2}\left(-i\pi\mathrm{sgn}\left(\mathrm{Im}\left(-z\right)\right)-\log\left(3\right)e^{-3z}+\mathrm{Ei}\left(-3z\right)\right)\\
+ & i\pi\mathrm{sgn}\left(\mathrm{Im}\left(-z\right)\right)+\log\left(2\right)e^{-2z}-\mathrm{Ei}\left(-2z\right)
\end{align*}
where $\gamma\left(\rho,2z\right)$ is the lower incomplete Gamma
function, $\rho=\beta+i\gamma$ runs over the non-trivial zeros of
the Riemann Zeta function and 
\[
\mathrm{Ei}\left(z\right)=-\int_{-z}^{\infty}\frac{e^{-t}}{t}dt,\ \left|\arg\left(z\right)\right|<\pi
\]
 is the Exponential integral function. Furthermore if we put
\[
\widetilde{S}\left(z\right)=\frac{e^{-2z}}{z}-\sum_{\rho}\Gamma\left(\rho\right)z^{-\rho}+\sum_{\rho}\left(z^{-\rho}\gamma\left(\rho,2z\right)-\frac{2^{\rho}e^{-2z}}{\rho}\right)+E\left(a,y\right)
\]
we get the bound
\[
E\left(a,y\right)\ll\begin{cases}
1+\left|z\right|^{1/2}, & \left|y\right|\leq a\\
1+\left|z\right|^{1/2}\left(1+\log^{2}\left(\left|y\right|/a\right)\right), & a<\left|y\right|\leq1\\
1, & \left|y\right|>1.
\end{cases}
\]
\end{thm}
Note that we improved the error term \eqref{errore s tilda} as $\left|y\right|\rightarrow\infty.$
This is interesting since, following the Languasco and Zaccagnini
approach, the size of the error term $E\left(a,y\right)$ is linked
to the size of $k$ (see Section $3$ of \cite{Langzac} for more
details). We also obtain the extra term $\sum_{\rho}\left(z^{-\rho}\gamma\left(\rho,2z\right)-\frac{2^{\rho}e^{-2z}}{\rho}\right)$;
following, again, the technique in \cite{Langzac}, we see that it
would be nice to have a sharp uniform estimation of this series. Unfortunately,
even if this series, as we will see, converges absolutely, it is not
simple to get a good uniform estimation. A deeper analysis of this
fact will be done in future research. Also note that the presence
of $e^{-2z}/z$ as ``main term'' instead of $1/z$ does not significantly
alter the asymptotic formula \eqref{eq:asint s tilda zac}. For example
if $z\in\mathbb{R}^{+}$ is clear that 
\[
\widetilde{S}\left(z\right)\sim e^{-2z}/z\sim1/z
\]
as $z\rightarrow0^{+},$ which is one of the form of the Prime Number
Theorem. In addition we have the following
\begin{thm}
\label{thm:k=00003D1}Let $N>4$ be a natural number. Then
\begin{align*}
\sum_{n\leq N}r_{G}\left(n\right)\left(N-n\right)= & \frac{N^{3}}{6}-2\sum_{\rho}\frac{\left(N-2\right)^{\rho+2}}{\rho\left(\rho+1\right)\left(\rho+2\right)}\\
+ & \sum_{\rho_{1}}\sum_{\rho_{2}}\frac{\Gamma\left(\rho_{1}\right)\Gamma\left(\rho_{2}\right)}{\Gamma\left(\rho_{1}+\rho_{2}+2\right)}N^{\rho_{1}+\rho_{2}+1}-2\sum_{\rho_{1}}\sum_{\rho_{2}}\frac{N^{\rho_{1}+\rho_{2}+1}}{\rho_{1}\rho_{2}}B_{2/N}\left(\rho_{1}+1,\rho_{2}+1\right)\\
+ & F\left(N\right)
\end{align*}
where $F\left(N\right)$ is a sum of (explicitly calculate) elementary
functions, dilogarithms and sum over non-trivial zeros of the Riemann
Zeta function involving the incomplete Beta function. Furthermore
we will prove that
\[
F\left(N\right)=O\left(N^{2}\right)
\]
as $N\rightarrow\infty.$
\end{thm}
This explicit formula extend, in some sense, the main result in \cite{Langzac}
(since their formula holds for $k>1$) and in \cite{GolYan} (since
the authors assume the Riemann Hypothesis). Furthermore it provides
different way to write the explicit formula of the Cesàro average
of Goldbach representations in the case $k=1$ respect to the main
formula in \cite{Perelli}. The study of the double series $\sum_{\rho_{1}}\sum_{\rho_{2}}\frac{N^{\rho_{1}+\rho_{2}+1}}{\rho_{1}\rho_{2}}B_{2/N}\left(\rho_{1}+1,\rho_{2}+1\right)$
is quite delicate since it is not obvious if converges absolutely
or not but it is clear that it plays a central role in formulae like
\eqref{eq:GoldeYan}. A deeper analysis of this type of double series
will be the subject of future research.

I thank my mentor Alessandro Zaccagnini, Jacopo ``Jack'' D'Aurizio
and Matthias Kunik for several conversation on this topic. 

\section{Notations}

To avoid ambiguity but keeping the standard notations, we specify
some symbols:

$\gamma,$ with or without subscripts, will be the imaginary part
of the non-trivial zeros of the Riemann Zeta function,

$\gamma\left(a,b\right)$ will be the lower incomplete Gamma function,

$\mathcal{L}\left(f\left(t\right)\right)\left(z\right)$ will be the
Laplace transform of $f\left(t\right).$

\section{Proof of Theorem 1}

We first recall that the upper incomplete Gamma function and the lower
incomplete Gamma function are defined as
\begin{equation}
\Gamma\left(a,z\right)=\int_{z}^{\infty}t^{a-1}e^{-t}dt\label{up gamma}
\end{equation}
\begin{equation}
\gamma\left(a,z\right)=\int_{0}^{z}t^{a-1}e^{-t}dt,\,\textrm{Re}\left(a\right)>0\label{low gamma}
\end{equation}
(see \cite{Olver}, formula $8.2.1$ and $8.2.2$ at page $174$).
Furthermore we have
\begin{equation}
\Gamma\left(a,z\right)+\gamma\left(a,z\right)=\Gamma\left(a\right)\label{eq:G1}
\end{equation}
\begin{equation}
\Gamma\left(a+1,z\right)=a\Gamma\left(a,z\right)+e^{-z}z^{a}\label{eq:G2}
\end{equation}
(see \cite{Olver}, formula $8.2.3$ at page $174$ and formula $8.8.2$
at page $178$). 

Let us consider the finite sum $\sum_{m=2}^{M}\Lambda\left(m\right)e^{-mz},\,M>2.$
By Abel summation we have
\[
\sum_{m=2}^{M}\Lambda\left(m\right)e^{-mz}=\psi\left(M\right)e^{-Mz}+z\int_{2}^{M}\psi\left(t\right)e^{-tz}dt
\]
where $\psi\left(M\right)$ is the Chebyshev psi function. From PNT
we have 
\[
\left|\psi\left(M\right)e^{-Mz}\right|\leq\psi\left(M\right)e^{-Ma}\ll Me^{-Ma}\rightarrow0
\]
as $M\rightarrow\infty$ so 
\begin{equation}
\widetilde{S}\left(z\right)=\sum_{m\geq1}\Lambda\left(m\right)e^{-mz}=z\int_{2}^{\infty}\psi\left(t\right)e^{-tz}dt\label{eq:Laplace psi}
\end{equation}
so $z^{-1}\widetilde{S}\left(z\right)$ is the Laplace transform of
$\psi\left(t\right).$ From the explicit formula for $\psi\left(t\right)$
(see for example \cite{daven}, page $104$) we obtain 
\begin{align}
z\int_{2}^{\infty}\psi\left(t\right)e^{-tz}dt= & z\int_{2}^{\infty}te^{-tz}dt-z\int_{2}^{\infty}\sum_{\rho}\frac{t^{\rho}}{\rho}e^{-tz}dt\nonumber \\
- & \frac{\zeta^{\prime}}{\zeta}\left(0\right)z\int_{2}^{\infty}e^{-tz}dt-\frac{z}{2}\int_{2}^{\infty}\log\left(1-t^{-2}\right)e^{-tz}dt.\label{eq:9}
\end{align}
Now we integrate term by term \eqref{eq:9} and we assume for now
that we can switch the series over the non-trivial zeros and the integral
(we will show that it is possible at the end of this proof). We have
immediately that
\[
z\int_{2}^{\infty}te^{-tz}dt=2e^{-2z}+\frac{e^{-2z}}{z}
\]
\[
-\frac{\zeta^{\prime}}{\zeta}\left(0\right)z\int_{2}^{\infty}e^{-tz}dt=-\frac{\zeta^{\prime}}{\zeta}\left(0\right)e^{-2z}.
\]
Now, from \eqref{eq:G1} and \eqref{eq:G2}, we have 
\begin{align}
-z\int_{2}^{\infty}\sum_{\rho}\frac{t^{\rho}}{\rho}e^{-tz}dt= & -\sum_{\rho}\frac{1}{\rho}z\int_{2}^{\infty}t^{\rho}e^{-tz}dt\nonumber \\
= & -\sum_{\rho}\frac{1}{\rho}z^{-\rho}\int_{2z}^{\infty}t^{\rho}e^{-t}dt\nonumber \\
= & -\sum_{\rho}\frac{1}{\rho}z^{-\rho}\Gamma\left(\rho+1,2z\right)\nonumber \\
= & -\sum_{\rho}\left(z^{-\rho}\Gamma\left(\rho,2z\right)+\frac{2^{\rho}e^{-2z}}{\rho}\right)\nonumber \\
= & -\sum_{\rho}\left(z^{-\rho}\Gamma\left(\rho\right)-z^{-\rho}\gamma\left(\rho,2z\right)+\frac{2^{\rho}e^{-2z}}{\rho}\right)\nonumber \\
= & -\sum_{\rho}z^{-\rho}\Gamma\left(\rho\right)+\sum_{\rho}\left(z^{-\rho}\gamma\left(\rho,2z\right)-\frac{2^{\rho}e^{-2z}}{\rho}\right).\label{eq:identit=0000E0}
\end{align}

Note that the series in \eqref{eq:identit=0000E0} converges absolutely;
for $-\sum_{\rho}z^{-\rho}\Gamma\left(\rho\right)$ see \cite{Kunik},
formula $(2.2)$ and for $\sum_{\rho}\left(z^{-\rho}\gamma\left(\rho,2z\right)-\frac{2^{\rho}e^{-2z}}{\rho}\right)$,
using the integral representation \eqref{low gamma} and integrating
by pars, we observe that
\[
\sum_{\rho}\left(z^{-\rho}\gamma\left(\rho,2z\right)-\frac{2^{\rho}e^{-2z}}{\rho}\right)=z\sum_{\rho}\frac{1}{\rho}\int_{0}^{2}t^{\rho}e^{-tz}dt\ll\left(\left|z\right|+\left|z\right|^{2}\right)\sum_{\rho}\frac{1}{\left|\rho\right|^{2}}.
\]
 Now we consider the integral involving $\log\left(1-t^{-2}\right).$
Integrating by parts, we observe that
\begin{equation}
-\frac{z}{2}\int\log\left(t\right)e^{-zt}dt=\frac{1}{2}\left(\log\left(t\right)e^{-zt}-\textrm{Ei}\left(-zt\right)\right)+K\label{eq:primitivaEi}
\end{equation}
where $\textrm{Ei}\left(x\right)$ is the Exponential integral function
and $K$ is a constant. So 
\begin{align*}
-\frac{z}{2}\int_{2}^{\infty}\log\left(1-t^{-2}\right)e^{-tz}dt= & -\frac{z}{2}\int_{2}^{\infty}\log\left(t-1\right)e^{-tz}dt-\frac{z}{2}\int_{2}^{\infty}\log\left(t+1\right)e^{-tz}dt\\
+ & z\int_{2}^{\infty}\log\left(t\right)e^{-tz}dt\\
= & \mathcal{J}_{1}+\mathcal{J}_{2}+\mathcal{J}_{3},
\end{align*}
say. Let us study $\mathcal{J}_{1}$. By \eqref{eq:primitivaEi} we
have that
\begin{align*}
\mathcal{J}_{1}= & -\frac{e^{-z}}{2}z\int_{1}^{\infty}\log\left(t\right)e^{-tz}dt\\
= & \frac{e^{-z}}{2}\left(\log\left(t\right)e^{-tz}-\textrm{Ei}\left(-tz\right)\right)_{1}^{\infty}.
\end{align*}
From the asymptotic
\begin{equation}
\textrm{Ei}\left(w\right)\sim i\pi\textrm{sgn}\left(\textrm{Im}\left(w\right)\right)-\frac{e^{w}}{w}\label{eq:eiasimpt}
\end{equation}
as $\left|w\right|\rightarrow\infty$ (see for example \cite{Peg},
equation $(26)$ at page $192$) we get 
\begin{equation}
\lim_{t\rightarrow\infty}\left(\log\left(t\right)e^{-tz}-\textrm{Ei}\left(-tz\right)\right)=-i\pi\textrm{sgn}\left(\textrm{Im}\left(-z\right)\right).\label{eq:limitEi}
\end{equation}
So
\[
\mathcal{J}_{1}=\frac{e^{-z}}{2}\left(-i\pi\textrm{sgn}\left(\textrm{Im}\left(-z\right)\right)+\textrm{Ei}\left(-z\right)\right).
\]
Let us consider $\mathcal{J}_{2}$. Using again \eqref{eq:limitEi}
we get
\begin{align*}
\mathcal{J}_{2}= & -\frac{e^{z}}{2}z\int_{3}^{\infty}\log\left(t\right)e^{-tz}dt\\
= & \frac{e^{z}}{2}\left(\log\left(t\right)e^{-tz}-\textrm{Ei}\left(-tz\right)\right)_{3}^{\infty}\\
= & \frac{e^{z}}{2}\left(-i\pi\textrm{sgn}\left(\textrm{Im}\left(-z\right)\right)-\log\left(3\right)e^{-3z}+\textrm{Ei}\left(-3z\right)\right).
\end{align*}
We have now to analyze $\mathcal{J}_{3}$. From \eqref{eq:limitEi}
we have
\begin{align*}
\mathcal{J}_{3}= & -\left(\log\left(t\right)e^{-tz}-\textrm{Ei}\left(-tz\right)\right)_{2}^{\infty}\\
= & i\pi\textrm{sgn}\left(\textrm{Im}\left(-z\right)\right)+\log\left(2\right)e^{-2z}-\textrm{Ei}\left(-2z\right).
\end{align*}
Finally
\begin{align*}
-\frac{z}{2}\int_{2}^{\infty}\log\left(1-t^{-2}\right)e^{-tz}dt= & \frac{e^{-z}}{2}\left(-i\pi\textrm{sgn}\left(\textrm{Im}\left(-z\right)\right)+\textrm{Ei}\left(-z\right)\right)\\
+ & \frac{e^{z}}{2}\left(-i\pi\textrm{sgn}\left(\textrm{Im}\left(-z\right)\right)-\log\left(3\right)e^{-3z}+\textrm{Ei}\left(-3z\right)\right)\\
+ & i\pi\textrm{sgn}\left(\textrm{Im}\left(-z\right)\right)+\log\left(2\right)e^{-2z}-\textrm{Ei}\left(-2z\right).
\end{align*}
So 
\begin{align*}
\widetilde{S}\left(z\right)= & 2e^{-2z}+\frac{e^{-2z}}{z}-\sum_{\rho}z^{-\rho}\Gamma\left(\rho\right)\\
+ & \sum_{\rho}\left(z^{-\rho}\gamma\left(\rho,2z\right)-\frac{2^{\rho}e^{-2z}}{\rho}\right)-\frac{\zeta^{\prime}}{\zeta}\left(0\right)e^{-2z}\\
+ & \frac{e^{-z}}{2}\left(-i\pi\textrm{sgn}\left(\textrm{Im}\left(-z\right)\right)+\textrm{Ei}\left(-z\right)\right)\\
+ & \frac{e^{z}}{2}\left(-i\pi\textrm{sgn}\left(\textrm{Im}\left(-z\right)\right)-\log\left(3\right)e^{-3z}+\textrm{Ei}\left(-3z\right)\right)\\
+ & i\pi\textrm{sgn}\left(\textrm{Im}\left(-z\right)\right)+\log\left(2\right)e^{-2z}-\textrm{Ei}\left(-2z\right).
\end{align*}
We have almost proved the first part of the theorem. Now we have to
show that we may exchange the series with the integral in \eqref{eq:9}.
Let us define
\begin{equation}
h_{T}\left(t\right)=e^{-tz}\sum_{\rho:\,\left|\gamma\right|\leq T}\frac{t^{\rho}}{\rho}\label{eq:convdominata}
\end{equation}
where $T>2.$ From the truncated explicit formula (see \cite{daven}
at page $109$) we have
\begin{equation}
\psi\left(t\right)=t-\sum_{\rho:\,\left|\gamma\right|\leq T}\frac{t^{\rho}}{\rho}+O\left(\frac{t\log^{2}\left(tT\right)}{T}+\log\left(t\right)\min\left(1,\frac{t}{T\left\langle t\right\rangle }\right)\right),\,t>2,\,T>2\label{psi troncata}
\end{equation}
where $\left\langle t\right\rangle $ is the distance from $t$ to
the nearest prime power, and so 
\begin{align*}
\left|h_{T}\left(t\right)\right|\ll & e^{-ta}\left(\left|\psi\left(t\right)-t\right|+\frac{t\log^{2}\left(tT\right)}{T}+\log\left(t\right)\right)\\
\ll & e^{-ta}\left(t+\frac{t\log^{2}\left(t\right)+2t\log\left(t\right)\log\left(T\right)+t\log^{2}\left(T\right)+T\log\left(t\right)}{T}\right)
\end{align*}
since $\left|\psi\left(t\right)-t\right|\ll t.$ But $\log^{2}\left(T\right)/T\ll1$
as $T\rightarrow\infty$, then 
\[
\left|h_{T}\left(t\right)\right|\ll e^{-ta}t\log^{2}\left(t\right)
\]
so we have to prove that 
\[
\int_{2}^{\infty}t\log^{2}\left(t\right)e^{-ta}dt
\]
converges. Trivially we have
\[
\int_{2}^{\infty}t\log^{2}\left(t\right)e^{-ta}dt\ll_{a}\int_{0}^{\infty}t^{2}e^{-t}dt=\Gamma\left(3\right)
\]
so by the dominated convergence theorem we may exchange the integral
with the series. It remains to prove the estimation of $E\left(a,y\right)$.
From Lemma 1 of \cite{Langzac} (The bound has been corrected in \cite{Lang})
we observe that
\[
\sum_{\rho}\left(z^{-\rho}\gamma\left(\rho,2z\right)-\frac{2^{\rho}e^{-2z}}{\rho}\right)\ll\begin{cases}
1+\left|z\right|^{1/2}, & \left|y\right|\leq a\\
1+\left|z\right|^{1/2}\left(1+\log^{2}\left(\left|y\right|/a\right)\right), & a<\left|y\right|\leq1
\end{cases}
\]
and
\[
E\left(a,y\right)\ll\begin{cases}
1+\left|z\right|^{1/2}, & \left|y\right|\leq a\\
1+\left|z\right|^{1/2}\left(1+\log^{2}\left(\left|y\right|/a\right)\right), & a<\left|y\right|\leq1
\end{cases}
\]
so the estimation follows. If $\left|y\right|>1$ from \eqref{eq:eiasimpt}
we see that 
\[
\left|\textrm{Ei}\left(-z\right)-i\pi\textrm{sgn}\left(\textrm{Im}\left(-z\right)\right)\right|\ll\left|\frac{e^{-z}}{z}\right|\ll1
\]
 as $y\rightarrow\pm\infty$ and this concludes the proof.

\section{Proof of Theorem 2}

Let $z=1/N+iy$, $N\in\mathbb{N}$, $N>4$ and $y\in\mathbb{R}.$
From \eqref{eq:Laplace psi} we have that
\[
\mathcal{L}\left(\psi\left(t\right)\right)\left(z\right)=z^{-1}\widetilde{S}\left(z\right)
\]
so by the convolution theorem and the inversion theorem of the Laplace
transform (see \cite{Folland} theorems $8.1$ and $8.4$), we get
\begin{align*}
\frac{1}{2\pi i}\int_{\left(1/N\right)}e^{uz}z^{-2}\widetilde{S}\left(z\right)^{2}dz= & \begin{cases}
\int_{2}^{u-2}\psi\left(t\right)\psi\left(u-t\right)dt, & u>4\\
0, & \mathrm{otherwise}
\end{cases}
\end{align*}
then from the explicit formula we have to evaluate
\[
\int_{2}^{u-2}\left(t-\sum_{\rho}\frac{t^{\rho}}{\rho}-\frac{\zeta^{\prime}}{\zeta}\left(0\right)-\frac{1}{2}\log\left(1-\frac{1}{t^{2}}\right)\right)\left(u-t-\sum_{\rho}\frac{\left(u-t\right)^{\rho}}{\rho}-\frac{\zeta^{\prime}}{\zeta}\left(0\right)-\frac{1}{2}\log\left(1-\frac{1}{\left(u-t\right)^{2}}\right)\right)dt
\]
\[
=\sum_{m=1}^{10}\int_{2}^{u-2}s_{m}\left(t,u\right)dt,
\]
say, where the functions $s_{m}\left(t,u\right)$ will be defined
in the next sections. 

\subsection{Integral of $s_{1}\left(t,u\right)$}

We take $s_{1}\left(t,u\right):=t\left(u-t\right).$ We get immediately
\[
\int_{2}^{u-2}s_{1}\left(t,u\right)dt=\frac{u^{3}-24u+32}{6}.
\]

\subsection{Integral of $s_{2}\left(t,u\right)$}

We take $s_{2}\left(t,u\right):=-2\left(u-t\right)\sum_{\rho}\frac{t^{\rho}}{\rho}.$
We have already seen in \eqref{psi troncata} that
\begin{equation}
\sum_{\rho:\,\left|\gamma\right|\leq T}\frac{t^{\rho}}{\rho}\ll t\log^{2}\left(t\right)\label{bound esplicito somma zeri troncata}
\end{equation}
then by the dominated convergence theorem we can write
\begin{align*}
\int_{2}^{u-2}s_{2}\left(t,u\right)dt= & -2\sum_{\rho}\frac{1}{\rho}\int_{2}^{u-2}\left(u-t\right)t^{\rho}dt\\
= & -4\sum_{\rho}\frac{\left(u-2\right)^{\rho+1}}{\rho\left(\rho+1\right)}+\left(u-2\right)\sum_{\rho}\frac{2^{\rho+2}}{\rho\left(\rho+1\right)}\\
- & 2\sum_{\rho}\frac{1}{\rho\left(\rho+1\right)}\int_{2}^{u-2}t^{\rho+1}dt\\
= & -4\sum_{\rho}\frac{\left(u-2\right)^{\rho+1}}{\rho\left(\rho+1\right)}+\left(u-2\right)\sum_{\rho}\frac{2^{\rho+2}}{\rho\left(\rho+1\right)}\\
- & 2\sum_{\rho}\frac{\left(u-2\right)^{\rho+2}}{\rho\left(\rho+1\right)\left(\rho+2\right)}+\sum_{\rho}\frac{2^{\rho+3}}{\rho\left(\rho+1\right)\left(\rho+2\right)}.
\end{align*}
Note that the all the series in this Section converges absolutely.

\subsection{Integral of $s_{3}\left(t,u\right)$}

We consider $s_{3}\left(t,u\right):=-2t\frac{\zeta^{\prime}}{\zeta}\left(0\right).$
We immediately obtain
\[
\int_{2}^{u-2}s_{3}\left(t,u\right)dt=\frac{\zeta^{\prime}}{\zeta}\left(0\right)\left(4u-u^{2}\right).
\]

\subsection{Integral of $s_{4}\left(t,u\right)$}

We take $s_{4}\left(t,u\right):=-\left(u-t\right)\log\left(1-\frac{1}{t^{2}}\right).$
We get
\begin{align*}
\int_{2}^{u-2}s_{4}\left(t,u\right)dt= & -\int_{2}^{u-2}\left(u-t\right)\log\left(1-\frac{1}{t^{2}}\right)dt\\
= & -\int_{2}^{u-2}\left(u-t\right)\log\left(t-1\right)dt-\int_{2}^{u-2}\left(u-t\right)\log\left(t+1\right)dt\\
+ & 2\int_{2}^{u-2}\left(u-t\right)\log\left(t\right)dt
\end{align*}
then, after simple but boring calculations, we obtain
\begin{align*}
\int_{2}^{u-2}s_{4}\left(t,u\right)dt= & \frac{1}{4}\left(\left(u-4\right)\left(3u-2\right)-2\left(u-3\right)\left(u+1\right)\log\left(u-3\right)\right.\\
- & 8+3u^{2}-6\log\left(3\right)+2u\left(\log\left(729\right)-5\right)-2\log\left(u-1\right)\left(u^{2}+2u-3\right)\\
- & \left.6u^{2}-8u\left(\log\left(4\right)-3\right)+2\log\left(256\right)+4\log\left(u-2\right)\left(u^{2}-4\right)\right).
\end{align*}
Also it is quite simple to observe that
\[
\int_{2}^{u-2}s_{4}\left(t,u\right)dt\ll u^{2}.
\]

\subsection{Integral of $s_{5}\left(t,u\right)$}

We now take $s_{5}\left(t,u\right):=\sum_{\rho_{1}}\frac{\left(u-t\right)^{\rho_{1}}}{\rho_{1}}\sum_{\rho_{2}}\frac{t^{\rho_{2}}}{\rho_{2}}.$
Then we have to calculate
\[
\int_{2}^{u-2}s_{5}\left(t,u\right)dt=\int_{2}^{u-2}\sum_{\rho_{1}}\frac{\left(u-t\right)^{\rho_{1}}}{\rho_{1}}\sum_{\rho_{2}}\frac{t^{\rho_{2}}}{\rho_{2}}dt.
\]
We have to prove that it is possible to exchange the integral and
the series over the non-trivial zeros of the Riemann Zeta function.
Let us define
\[
S_{T,u}^{1}\left(t\right)=\sum_{\rho_{1}:\,\left|\gamma_{1}\right|\leq T}\frac{\left(u-t\right)^{\rho_{1}}}{\rho_{1}}\sum_{\rho_{2}}\frac{t^{\rho_{2}}}{\rho_{2}},\,T>2.
\]
From \eqref{bound esplicito somma zeri troncata} and the trivial
bound
\begin{equation}
\left|\sum_{\rho_{2}}\frac{t^{\rho_{2}}}{\rho_{2}}\right|=\left|t-\psi\left(t\right)-\frac{\zeta^{\prime}}{\zeta}\left(0\right)-\frac{1}{2}\log\left(1-\frac{1}{t^{2}}\right)\right|\ll t\label{trivial}
\end{equation}
for $t>2,$ we have
\[
\left|S_{T,u}^{1}\left(t\right)\right|\ll t\left(u-t\right)\log^{2}\left(u-t\right)
\]
then by the dominated convergence theorem
\[
\int_{2}^{u-2}\sum_{\rho_{1}}\frac{\left(u-t\right)^{\rho_{1}}}{\rho_{1}}\sum_{\rho_{2}}\frac{t^{\rho_{2}}}{\rho_{2}}dt=\sum_{\rho_{1}}\frac{1}{\rho_{1}}\int_{2}^{u-2}\left(u-t\right)^{\rho_{1}}\sum_{\rho_{2}}\frac{t^{\rho_{2}}}{\rho_{2}}dt.
\]
Now if we take
\[
S_{T,u,\rho_{1}}^{2}\left(t\right)=\left(u-t\right)^{\rho_{1}}\sum_{\rho_{2}:\,\left|\gamma_{2}\right|\leq T}\frac{t^{\rho_{2}}}{\rho_{2}},\,T>2
\]
from \eqref{bound esplicito somma zeri troncata} we have
\[
\left|S_{T,u,\rho_{1}}^{2}\left(t\right)\right|\ll\left(u-t\right)t\log^{2}\left(t\right)
\]
then again by the dominated convergence theorem we can conclude that
\[
\sum_{\rho_{1}}\frac{1}{\rho_{1}}\int_{2}^{u-2}\left(u-t\right)^{\rho_{1}}\sum_{\rho_{2}}\frac{t^{\rho_{2}}}{\rho_{2}}dt=\sum_{\rho_{1}}\sum_{\rho_{2}}\frac{1}{\rho_{1}\rho_{2}}\int_{2}^{u-2}\left(u-t\right)^{\rho_{1}}t^{\rho_{2}}dt.
\]

Then
\begin{align*}
\int_{2}^{u-2}\sum_{\rho_{1}}\frac{\left(u-t\right)^{\rho_{1}}}{\rho_{1}}\sum_{\rho_{2}}\frac{t^{\rho_{2}}}{\rho_{2}}dt= & \sum_{\rho_{1}}\sum_{\rho_{2}}\frac{u^{\rho_{1}+\rho_{2}+1}}{\rho_{1}\rho_{2}}\frac{\Gamma\left(\rho_{1}+1\right)\Gamma\left(\rho_{2}+1\right)}{\Gamma\left(\rho_{1}+\rho_{2}+2\right)}\\
- & \sum_{\rho_{1}}\sum_{\rho_{2}}\frac{u^{\rho_{1}+\rho_{2}+1}}{\rho_{1}\rho_{2}}B_{2/u}\left(\rho_{1}+1,\rho_{2}+1\right)\\
- & \sum_{\rho_{1}}\sum_{\rho_{2}}\frac{u^{\rho_{1}+\rho_{2}+1}}{\rho_{1}\rho_{2}}B_{2/u}\left(\rho_{2}+1,\rho_{1}+1\right)\\
= & \sum_{\rho_{1}}\sum_{\rho_{2}}u^{\rho_{1}+\rho_{2}+1}\frac{\Gamma\left(\rho_{1}\right)\Gamma\left(\rho_{2}\right)}{\Gamma\left(\rho_{1}+\rho_{2}+2\right)}-2\sum_{\rho_{1}}\sum_{\rho_{2}}\frac{u^{\rho_{1}+\rho_{2}+1}}{\rho_{1}\rho_{2}}B_{2/u}\left(\rho_{1}+1,\rho_{2}+1\right)
\end{align*}
where $B_{x}\left(a,b\right)$ is the incomplete Beta function (for
more details see for example \cite{Olver}, ch. $8.17$). We recall
that 
\[
\sum_{\rho_{1}}\sum_{\rho_{2}}u^{\rho_{1}+\rho_{2}+1}\frac{\Gamma\left(\rho_{1}\right)\Gamma\left(\rho_{2}\right)}{\Gamma\left(\rho_{1}+\rho_{2}+2\right)}
\]
converges absolutely (see \cite{Langzac}), so $\sum_{\rho_{1}}\sum_{\rho_{2}}\frac{u^{\rho_{1}+\rho_{2}+1}}{\rho_{1}\rho_{2}}B_{2/u}\left(\rho_{1}+1,\rho_{2}+1\right)$
converges (at least conditionally), since
\begin{align*}
\left|2\sum_{\rho_{1}}\sum_{\rho_{2}}\frac{u^{\rho_{1}+\rho_{2}+1}}{\rho_{1}\rho_{2}}B_{2/u}\left(\rho_{1}+1,\rho_{2}+1\right)\right|\ll & \sum_{\rho_{1}}\sum_{\rho_{2}}\left|u^{\rho_{1}+\rho_{2}+1}\frac{\Gamma\left(\rho_{1}\right)\Gamma\left(\rho_{2}\right)}{\Gamma\left(\rho_{1}+\rho_{2}+2\right)}\right|\\
+ & \left|\int_{2}^{u-2}\sum_{\rho_{1}}\frac{\left(u-t\right)^{\rho_{1}}}{\rho_{1}}\sum_{\rho_{2}}\frac{t^{\rho_{2}}}{\rho_{2}}dt\right|
\end{align*}
and trivially
\[
\left|\int_{2}^{u-2}\sum_{\rho_{1}}\frac{\left(u-t\right)^{\rho_{1}}}{\rho_{1}}\sum_{\rho_{2}}\frac{t^{\rho_{2}}}{\rho_{2}}dt\right|\ll\int_{2}^{u-2}\left(u-t\right)tdt<\infty
\]
from \eqref{trivial}.

\subsection{Integral of $s_{6}\left(t,u\right)$}

We take $s_{6}\left(t,u\right):=2\frac{\zeta^{\prime}}{\zeta}\left(0\right)\sum_{\rho}\frac{t^{\rho}}{\rho}.$
From \eqref{bound esplicito somma zeri troncata} we know that it
is possible to switch the integral and the series, so
\begin{align*}
\int_{2}^{u-2}s_{6}\left(t,u\right)dt= & 2\frac{\zeta^{\prime}}{\zeta}\left(0\right)\sum_{\rho}\frac{1}{\rho}\int_{2}^{u-2}t^{\rho}dt\\
= & 2\frac{\zeta^{\prime}}{\zeta}\left(0\right)\sum_{\rho}\frac{\left(u-2\right)^{\rho+1}}{\rho\left(\rho+1\right)}-\frac{\zeta^{\prime}}{\zeta}\left(0\right)\sum_{\rho}\frac{2^{\rho+2}}{\rho\left(\rho+1\right)}.
\end{align*}
Again it is quite simple to observe that
\[
\int_{2}^{u-2}s_{6}\left(t,u\right)dt\ll u^{2}
\]
due to the absolute convergence of $\sum_{\rho}\frac{1}{\rho\left(\rho+1\right)}.$

\subsection{Integral of $s_{7}\left(t,u\right)$}

We consider $s_{7}\left(t,u\right):=\sum_{\rho}\frac{\left(u-t\right)^{\rho}}{\rho}\log\left(1-\frac{1}{t^{2}}\right)$.
Using the dominated convergence theorem we can prove that we can exchange
the integral and the series, so
\begin{align*}
\int_{2}^{u-2}s_{7}\left(t,u\right)dt= & \sum_{\rho}\frac{1}{\rho}\int_{2}^{u-2}\left(u-t\right)^{\rho}\log\left(1-\frac{1}{t^{2}}\right)dt\\
= & \sum_{\rho}\frac{1}{\rho}\int_{2}^{u-2}\left(u-t\right)^{\rho}\log\left(t-1\right)dt+\sum_{\rho}\frac{1}{\rho}\int_{2}^{u-2}\left(u-t\right)^{\rho}\log\left(t+1\right)dt\\
- & 2\sum_{\rho}\frac{1}{\rho}\int_{2}^{u-2}\left(u-t\right)^{\rho}\log\left(t\right)dt.
\end{align*}
Let us consider the first integral. Integrating by parts we get
\begin{align*}
H_{1}\left(u\right):= & \sum_{\rho}\frac{1}{\rho}\int_{2}^{u-2}\left(u-t\right)^{\rho}\log\left(t-1\right)dt\\
+ & \sum_{\rho}\frac{2^{\rho+1}\log\left(u-3\right)}{\rho\left(\rho+1\right)}-\sum_{\rho}\frac{1}{\rho\left(\rho+1\right)}\int_{2}^{u-2}\left(u-t\right)^{\rho+1}\left(t-1\right)^{-1}dt\\
= & \sum_{\rho}\frac{2^{\rho+1}\log\left(u-3\right)}{\rho\left(\rho+1\right)}-\sum_{\rho}\frac{\left(u-1\right)^{\rho+1}}{\rho\left(\rho+1\right)}\int_{2/\left(u-1\right)}^{1-1/\left(u-1\right)}s^{\rho+1}\left(1-s\right)^{-1}ds\\
= & \sum_{\rho}\frac{2^{\rho+1}\log\left(u-3\right)}{\rho\left(\rho+1\right)}-\sum_{\rho}\frac{\left(u-1\right)^{\rho+1}}{\rho\left(\rho+1\right)}\left[B_{1-1/\left(u-1\right)}\left(\rho+2,0\right)-B_{2/\left(u-1\right)}\left(\rho+2,0\right)\right]
\end{align*}
 In a similar way we have
\begin{align*}
H_{2}\left(u\right):= & \sum_{\rho}\frac{1}{\rho}\int_{2}^{u-2}\left(u-t\right)^{\rho}\log\left(t+1\right)dt\\
= & \sum_{\rho}\frac{2^{\rho+1}\log\left(u-1\right)}{\rho\left(\rho+1\right)}-\sum_{\rho}\frac{\left(u-2\right)^{\rho+1}\log\left(3\right)}{\rho\left(\rho+1\right)}\\
- & \sum_{\rho}\frac{1}{\rho\left(\rho+1\right)}\int_{2}^{u-2}\left(u-t\right)^{\rho+1}\left(t+1\right)^{-1}dt\\
= & \sum_{\rho}\frac{2^{\rho}\log\left(u-1\right)}{\rho\left(\rho+1\right)}-\sum_{\rho}\frac{\left(u-2\right)^{\rho}\log\left(3\right)}{\rho\left(\rho+1\right)}\\
- & \sum_{\rho}\frac{\left(u+1\right)^{\rho+1}}{\rho\left(\rho+1\right)}\left[B_{1-3/\left(u+1\right)}\left(\rho+2,0\right)-B_{2/\left(u+1\right)}\left(\rho+2,0\right)\right]
\end{align*}
and
\begin{align*}
H_{3}\left(u\right):= & -2\sum_{\rho}\frac{1}{\rho}\int_{2}^{u-2}\left(u-t\right)^{\rho}\log\left(t\right)dt\\
= & -\sum_{\rho}\frac{2^{\rho+2}\log\left(u-2\right)}{\rho\left(\rho+1\right)}+2\sum_{\rho}\frac{\left(u-2\right)^{\rho+1}\log\left(2\right)}{\rho\left(\rho+1\right)}\\
+ & 2\sum_{\rho}\frac{1}{\rho\left(\rho+1\right)}\int_{2}^{u-2}\frac{\left(u-t\right)^{\rho+1}}{t}dt\\
= & -\sum_{\rho}\frac{2^{\rho+2}\log\left(u-2\right)}{\rho\left(\rho+1\right)}+2\sum_{\rho}\frac{\left(u-2\right)^{\rho+1}\log\left(2\right)}{\rho\left(\rho+1\right)}\\
+ & 2\sum_{\rho}\frac{u^{\rho+1}}{\rho\left(\rho+1\right)}\left[B_{1-2/u}\left(\rho+2,0\right)-B_{2/u}\left(\rho+2,0\right)\right]
\end{align*}
Finally we can write
\[
\int_{2}^{u-2}s_{7}\left(t,u\right)dt=H_{1}\left(u\right)+H_{2}\left(u\right)+H_{3}\left(u\right).
\]
From the well known estimation
\[
\sum_{\rho}\frac{x^{\rho}}{\rho}\ll x\exp\left(-C\sqrt{\log\left(x\right)}\right),\,x>1,\,C>0
\]
(see \cite{daven}, chapter $18$) we can easily conclude that
\[
\int_{2}^{u-2}s_{7}\left(t,u\right)dt\ll u^{2}.
\]
Again we can observe that all the series over the non-trivial zeros
of $\zeta\left(s\right)$ in this Section converges absolutely; to
see this fact in the case of the series involving the incomplete Beta
function (which is, probably, the less evident absolute convergence)
note that
\[
2\sum_{\rho}\frac{u^{\rho+1}}{\rho\left(\rho+1\right)}\left[B_{1-2/u}\left(\rho+2,0\right)-B_{2/u}\left(\rho+2,0\right)\right]=2\sum_{\rho}\frac{1}{\rho\left(\rho+1\right)}\int_{2}^{u-2}\frac{\left(u-t\right)^{\rho+1}}{t}dt
\]
\[
\ll\sum_{\rho}\frac{1}{\left|\rho\right|\left|\rho+1\right|}\int_{2}^{u-2}\frac{\left(u-t\right)^{2}}{t}dt<\infty.
\]

\subsection{Integral of $s_{8}\left(t,u\right)$}

We take $s_{8}\left(t,u\right):=\frac{\zeta^{\prime}}{\zeta}\left(0\right)^{2}.$
We immediately get
\[
\int_{2}^{u-2}s_{8}\left(t,u\right)dt=\frac{\zeta^{\prime}}{\zeta}\left(0\right)^{2}\left(u-4\right).
\]

\subsection{Integral of $s_{9}\left(t,u\right)$}

We consider $s_{9}\left(t,u\right):=\frac{\zeta^{\prime}}{\zeta}\left(0\right)\log\left(1-\frac{1}{t^{2}}\right).$
We obtain
\begin{align*}
\int_{2}^{u-2}s_{9}\left(t,u\right)dt= & \frac{\zeta^{\prime}}{\zeta}\left(0\right)\left(\left(u-3\right)\log\left(u-3\right)\right.\\
+ & \log\left(\frac{1}{27\left(u-1\right)}\right)+u\log\left(u-1\right)\\
+ & \left.2\log\left(4\right)-2\left(u-2\right)\log\left(u-2\right)\right)
\end{align*}
and again we trivially get the bound
\[
\int_{2}^{u-2}s_{9}\left(t,u\right)dt\ll u^{2}.
\]

\subsection{Integral of $s_{10}\left(t,u\right)$}

We take $s_{10}\left(t,u\right):=\frac{1}{4}\log\left(1-\frac{1}{t^{2}}\right)\log\left(1-\frac{1}{\left(u-t\right)^{2}}\right).$
We have that
\begin{align}
\int_{2}^{u-2}s_{10}\left(t,u\right)dt=\frac{1}{4} & \int_{2}^{u-2}\log\left(1-\frac{1}{t^{2}}\right)\log\left(1-\frac{1}{\left(u-t\right)^{2}}\right)dt\nonumber \\
= & \frac{1}{4}\int_{2}^{u-2}\log\left(\frac{t^{2}-1}{t^{2}}\right)\log\left(\frac{\left(u-t\right)^{2}-1}{\left(u-t\right)^{2}}\right)dt\nonumber \\
= & \frac{1}{4}\int_{2}^{u-2}\left(\log\left(t-1\right)+\log\left(t+1\right)-2\log\left(t\right)\right)\nonumber \\
\times & \left(\log\left(u-t-1\right)+\log\left(u-t+1\right)-2\log\left(u-t\right)\right)dt\label{integrali prodotti log}
\end{align}
so we have to evaluate all the combinations of \eqref{integrali prodotti log}.
We will calculate explicitly only the first, since the others are
similar. We consider 
\begin{align*}
\int_{2}^{u-2}\log\left(t-1\right)\log\left(u-t-1\right)dt= & \int_{1}^{u-3}\log\left(v\right)\log\left(u-v-2\right)dv\\
= & \left(u-2\right)\int_{1/\left(u-2\right)}^{1-1/\left(u-2\right)}\log\left(\left(u-2\right)s\right)\log\left(\left(u-2\right)\left(1-s\right)\right)ds\\
= & \left(u-2\right)\log^{2}\left(u-2\right)\left(1-\frac{2}{u-2}\right)\\
+ & \left(u-2\right)\log\left(u-2\right)\int_{1/\left(u-2\right)}^{1-1/\left(u-2\right)}\log\left(1-s\right)ds\\
+ & \left(u-2\right)\log\left(u-2\right)\int_{1/\left(u-2\right)}^{1-1/\left(u-2\right)}\log\left(s\right)ds\\
+ & \left(u-2\right)\int_{1/\left(u-2\right)}^{1-1/\left(u-2\right)}\log\left(s\right)\log\left(1-s\right)ds.
\end{align*}
Maybe it is useful to recall that, integrating by parts, we get
\begin{align*}
\left(u-2\right)\int_{1/\left(u-2\right)}^{1-1/\left(u-2\right)}\log\left(s\right)\log\left(1-s\right)ds= & \left(u-4\right)\log\left(1-\frac{1}{u-2}\right)\log\left(\frac{1}{u-2}\right)\\
- & \left(u-2\right)\int_{1/\left(u-2\right)}^{1-1/\left(u-2\right)}s\left(\frac{\log\left(1-s\right)}{s}-\frac{\log\left(s\right)}{1-s}\right)ds
\end{align*}
and
\begin{align*}
-\left(u-2\right)\int_{1/\left(u-2\right)}^{1-1/\left(u-2\right)}s\left(\frac{\log\left(1-s\right)}{s}-\frac{\log\left(s\right)}{1-s}\right)ds= & -\left(u-2\right)\int_{1/\left(u-2\right)}^{1-1/\left(u-2\right)}\log\left(1-s\right)ds\\
- & \left(u-2\right)\int_{1/\left(u-2\right)}^{1-1/\left(u-2\right)}\log\left(s\right)ds\\
+ & \left(u-2\right)\int_{1/\left(u-2\right)}^{1-1/\left(u-2\right)}\frac{\log\left(s\right)}{1-s}ds\\
= & -2\left(u-3\right)\log\left(u-3\right)+2\left(u-4\right)\left(1+\log\left(u-2\right)\right)\\
+ & \left(u-2\right)\left(\mathrm{Li}_{2}\left(\frac{1}{u-2}\right)-\mathrm{Li}_{2}\left(1-\frac{1}{u-2}\right)\right)
\end{align*}
where $\mathrm{Li}_{2}\left(x\right)$ is the Dilogartihm function.
After a boring calculation we obtain
\begin{align*}
V_{1}\left(u\right):=\frac{1}{4}\int_{2}^{u-2}\log\left(t-1\right)\log\left(u-t-1\right)dt= & \frac{1}{4}\left(2\left(u-4\right)+\log\left(u-3\right)\left(6-2u+\left(u-2\right)\log\left(u-2\right)\right)\right.\\
+ & \left.\left(u-2\right)\mathrm{Li}_{2}\left(\frac{1}{u-2}\right)-\left(u-2\right)\mathrm{Li}_{2}\left(1-\frac{1}{u-2}\right)\right).
\end{align*}
For the explicit calculations of the other integrals see the appendix.
Furthermore it is clear that
\[
\int_{2}^{u-2}s_{10}\left(t,u\right)dt\ll u^{2}
\]
holds. 

To finish the proof we have only to take $u=N$, recalling that 
\[
\sum_{n\leq N}r_{G}\left(n\right)\left(N-n\right)=\frac{1}{2\pi i}\int_{\left(1/N\right)}e^{Nz}z^{-2}\widetilde{S}\left(z\right)^{2}dz.
\]

\section{Appendix}

In this section we calculate explicitly all the integrals from the
Section $4.10$. Using the same ideas that we used for $V_{1}\left(u\right)$
(and the Mathematica's help) we get
\begin{align*}
V_{2}\left(u\right):= & \frac{1}{4}\int_{2}^{u-2}\log\left(t-1\right)\log\left(u-t+1\right)dt\\
= & \frac{1}{4}\left(-8+2u+\log\left(27\right)+\log\left(u-1\right)+\log\left(u-3\right)\left(u\log\left(u\right)+3-3\log\left(3\right)\right)\right.\\
- & \left.u\log\left(u^{2}-4u+3\right)+u\mathrm{Li}_{2}\left(\frac{1}{u}\right)-u\mathrm{Li}_{2}\left(1-\frac{3}{u}\right)\right),
\end{align*}
\begin{align*}
V_{3}\left(u\right):= & \frac{1}{4}\int_{2}^{u-2}\log\left(t+1\right)\log\left(u-t-1\right)dt\\
= & \frac{1}{4}\left(-8+2u-3\left(-1+\log\left(3\right)\right)\log\left(u-3\right)+\log\left(27(u-1)\right)+u\log\left(3\right)\log\left(1-\frac{3}{u}\right)\right.\\
+ & \left.u\log\left(u-1\right)\log\left(u\right)-u\log\left(u^{2}-4u+3\right)+u\mathrm{Li}_{2}\left(\frac{3}{u}\right)-u\mathrm{Li}_{2}\left(1-\frac{1}{u}\right)\right),
\end{align*}
\begin{align*}
V_{4}\left(u\right):= & \frac{1}{4}\int_{2}^{u-2}\log\left(t+1\right)\log\left(u-t+1\right)dt\\
= & \frac{1}{4}\left(-8+2u+\log\left(729\right)+2\log\left(u-1\right)-2u\log\left(u-1\right)-4\log\left(3\right)\log\left(u-1\right)\right.\\
+ & u\log\left(3\right)\log\left(u-1\right)-u\log\left(3\right)\log\left(u+2\right)-\log\left(9\right)\log\left(u+2\right)+2\log\left(u-1\right)\log\left(u+2\right)\\
- & \left.u\log\left(u-1\right)\log\left(u+2\right)+\left(2+u\right)\mathrm{Li}_{2}\left(\frac{3}{u+2}\right)-\left(2+u\right)\mathrm{Li}_{2}\left(1-\frac{3}{u+2}\right)\right),
\end{align*}
\begin{align*}
V_{5}\left(u\right):= & -\frac{1}{2}\int_{2}^{u-2}\log\left(t-1\right)\log\left(u-t\right)dt\\
= & \frac{1}{2}\left(8-2u-\log\left(4\right)-2\log\left(u-2\right)-\log\left(u-3\right)\left(3+u\log\left(u-1\right)-\log\left(4\left(u-1\right)\right)\right)\right.\\
+ & \left.u\log\left(u^{2}-5u+6\right)-\left(u-1\right)\mathrm{Li}_{2}\left(\frac{1}{u-1}\right)+\left(u-1\right)\mathrm{Li}_{2}\left(1-\frac{2}{u-1}\right)\right)
\end{align*}
\begin{align*}
V_{6}\left(u\right):= & -\frac{1}{2}\int_{2}^{u-2}\log\left(t+1\right)\log\left(u-t\right)dt\\
= & \frac{1}{2}\left(8-2u-2\log\left(u-2\right)+u\log\left(u-2\right)-u\log\left(3\right)\log\left(u-2\right)\right.\\
+ & \log\left(9\right)\log\left(u-2\right)+\log\left(\frac{1}{108\left(u-1\right)}\right)+u\log\left(u-1\right)+\log\left(4\right)\log\left(u-1\right)\\
+ & \log\left(3\right)\log\left(u+1\right)+u\log\left(3\right)\log\left(u+1\right)-\log\left(u-1\right)\log\left(u+1\right)\\
- & \left.u\log\left(u-1\right)\log\left(u+1\right)-\left(u+1\right)\mathrm{Li}_{2}\left(\frac{3}{u+1}\right)+\left(u+1\right)\mathrm{Li}_{2}\left(1-\frac{2}{u+1}\right)\right),
\end{align*}
\begin{align*}
V_{7}\left(u\right):= & -\frac{1}{2}\int_{2}^{u-2}\log\left(t\right)\log\left(u-t-1\right)dt\\
= & \frac{1}{2}\left(8-2u-\log\left(4\right)-3\log\left(u-3\right)+u\log\left(u-3\right)-u\log\left(2\right)\log\left(u-3\right)\right.\\
+ & \log\left(8\right)\log\left(u-3\right)-2\log\left(u-2\right)+u\log\left(u-2\right)-\log\left(2\right)\log\left(u-1\right)\\
+ & u\log\left(2\right)\log\left(u-1\right)+\log\left(u-2\right)\log\left(u-1\right)-u\log\left(u-2\right)\log\left(u-1\right)\\
- & \left.\left(u-1\right)\mathrm{Li}_{2}\left(\frac{2}{u-1}\right)+\left(u-1\right)\mathrm{Li}_{2}\left(1-\frac{1}{u-1}\right)\right),
\end{align*}
\begin{align*}
V_{8}\left(u\right):= & -\frac{1}{2}\int_{2}^{u-2}\log\left(t\right)\log\left(u-t+1\right)dt\\
= & \frac{1}{2}\left(8-2u-2\log\left(u-2\right)+u\log\left(u-2\right)]+\log\left(27\right)\log\left(u-2\right)+\log\left(\frac{1}{108\left(u-1\right)}\right)\right.\\
+ & u\log\left(u-1\right)-\log\left(2\right)\log\left(u-1\right)-u\log\left(2\right)\log\left(u-1\right)+\log\left(4\right)\log\left(u-1\right)\\
+ & \log\left(2\right)\log\left(u+1\right)+u\log\left(2\right)\log\left(u+1\right)-\log\left(u-2\right)\log\left(u+1\right)\\
- & \left.u\log\left(u-2\right)\log\left(u+1\right)-\left(u+1\right)\mathrm{Li}_{2}\left(\frac{2}{u+1}\right)+\left(u+1\right)\mathrm{Li}_{2}\left(1-\frac{3}{u+1}\right)\right)
\end{align*}
and finally
\begin{align*}
V_{9}\left(u\right):=\int_{2}^{u-2}\log\left(t\right)\log\left(u-t\right)dt= & -8+2u+\log\left(16\right)-u\log\left(2\right)\log\left(u\right)\\
+ & \log\left(u-2\right)\left(-2\left(-2+u+\log\left(4\right)\right)+u\log\left(2u\right)\right)\\
+ & u\left(\mathrm{Li}_{2}\left(\frac{2}{u}\right)-\mathrm{Li}_{2}\left(1-\frac{2}{u}\right)\right).
\end{align*}
So we can write
\[
\int_{2}^{u-2}s_{10}\left(t,u\right)dt=\sum_{m=1}^{9}V_{m}\left(u\right).
\]

\end{document}